\documentclass{amsart}
\usepackage{graphicx}

\usepackage[centertags]{amsmath}
\usepackage{amsfonts}
\usepackage{amssymb}
\usepackage{amsthm}
\usepackage{newlfont}

\newtheorem{thm}{Theorem}[section]
\newtheorem{cor}[thm]{Corollary}
\newtheorem{lem}[thm]{Lemma}
\newtheorem{prop}[thm]{Proposition}

\theoremstyle{definition}
\newtheorem{defn}[thm]{Definition}
\newtheorem{prob}[thm]{Problem}
\theoremstyle{remark}

\numberwithin{equation}{section}

\newcommand{\N}{\mathbb{N}}
\newcommand{\Z}{\mathbb{Z}}
\newcommand{\F}{\mathbb{F}}
\newcommand{\V}{\mathrm{V}}
\newcommand{\E}{\mathrm{E}}

\newcommand{\Cay}{\mathrm{Cay}}
\newcommand{\dom}{\mathrm{dom}}

\newcommand{\id}{\mathrm{id}}
\newcommand{\rmN}{\mathcal{N}}

\renewcommand{\:}{\,:\,}
\def\co{\colon\thinspace}

\begin{document}

\title{Burnside's Problem, Spanning Trees, and Tilings}
\author{Brandon Seward}
\address{Department of Mathematics, University of Michigan, 530 Church Street, Ann Arbor, MI 48109, U.S.A.}
\email{b.m.seward@gmail.com}

\begin{abstract}
In this paper we study geometric versions of Burnside's Problem and the von Neumann Conjecture. This is done by considering the notion of a translation-like action. Translation-like actions were introduced by Kevin Whyte as a geometric analogue of subgroup containment. Whyte proved a geometric version of the von Neumann Conjecture by showing that a finitely generated group is non-amenable if and only if it admits a translation-like action by any (equivalently every) non-abelian free group. We strengthen Whyte's result by proving that this translation-like action can be chosen to be transitive when the acting free group is finitely generated. We furthermore prove that the geometric version of Burnside's Problem holds true. That is, every finitely generated infinite group admits a translation-like action by $\Z$. This answers a question posed by Whyte. In pursuit of these results we discover an interesting property of Cayley graphs: every finitely generated infinite group $G$ has some Cayley graph having a regular spanning tree. This regular spanning tree can be chosen to have degree $2$ (and hence be a bi-infinite Hamiltonian path) if and only if $G$ has finitely many ends, and it can be chosen to have any degree greater than $2$ if and only if $G$ is non-amenable. We use this last result to then study tilings of groups. We define a general notion of polytilings and extend the notion of MT groups and ccc groups to the setting of polytilings. We prove that every countable group is poly-MT and every finitely generated group is poly-ccc.
\end{abstract}
\maketitle

\section{Introduction} \label{SEC INTRO}

This paper focuses on the notion of a translation-like action, introduced by Kevin Whyte in \cite{KW}, and its relevance to geometric versions of Burnside's Problem and the von Neumann Conjecture. The classical Burnside's Problem asks if every finitely generated infinite group contains $\Z$ as a subgroup, and the classical von Neumann Conjecture states that a group is non-amenable if and only if it contains a non-abelian free group as a subgroup. Although both of these problems are known to have negative answers (by work of Golod--Shafarevich \cite{GS} and Olshanskii \cite{O}, respectively), translation-like actions provide us with a new geometric perspective of these classical problems. Study of this notion in turn leads us to findings about Hamiltonian paths, regular spanning trees of Cayley graphs, and tilings of groups. While we are primarily interested in translation-like actions on groups, we are naturally led to consider translation-like actions on graphs as well. We therefore present a general definition of translation-like actions.

\begin{defn}[Whyte, \cite{KW}] \label{DEF TRANS ACT}
Let $H$ be a group and let $(X, d)$ be a metric space. A right action $*$ of $H$ on $X$ is \emph{translation-like} if it satisfies the following two conditions:
\begin{enumerate}
\item[\rm (i)] the action is free (i.e. $x * h = x$ implies $h = 1_H$);
\item[\rm (ii)] for every $h \in H$ the set $\{d(x, x * h) \: x \in X\}$ is bounded;
\end{enumerate}
\end{defn}

We view every finitely generated group as a metric space by using a left-invariant word-length metric associated to some finite generating set, and we view every connected graph as a metric space by using the path-length metric (see the next section for details).

We point out that if $H$ acts translation-like on a finitely generated group or a graph $X$, then for all $h \in H$ the map $x \in X \mapsto x * h \in X$ is bilipschitz \cite{KW}. Thus, in the context of finitely generated groups and graphs one could equivalently define translation-like actions to be free actions by bilipschitz maps at bounded distance from the identity map.

Whyte's motivation in defining translation-like actions is that it serves as a geometric analogue of subgroup containment. Specifically, we have the following proposition whose proof is trivial.

\begin{prop}[Whyte, \cite{KW}] \label{SUBGROUP ACTION}
Let $G$ be a finitely generated group and let $H \leq G$. Then the natural right action of $H$ on $G$ is a translation-like action.
\end{prop}

In \cite{KW}, Whyte suggested that conjectures relating geometric properties to subgroup containment may hold true if one requires a translation-like action in place of a subgroup. Specifically, Whyte mentioned the following three well known conjectures/problems.

\begin{enumerate}
\item[(1)] Burnside's Problem: does every finitely generated infinite group contain $\Z$ as a subgroup?
\item[(2)] The von Neumann Conjecture: a group is non-amenable if and only if it contains a non-abelian free subgroup.
\item[(3)] Gersten Conjecture: a finitely generated group is not word hyperbolic if and only if it contains some Baumslag-Solitar group.
\end{enumerate}

All three of the above conjectures/problems were answered negatively by Golod--Shafarevich \cite{GS}, Olshanskii \cite{O}, and Brady \cite{NB}, respectively. Whyte suggested the following ``geometric reformulations'' of these problems.

\begin{enumerate}
\item[(1')] Geometric Burnside's Problem: does every finitely generated infinite group admit a translation-like action by $\Z$?
\item[(2')] Geometric von Neumann Conjecture: a finitely generated group is non-amenable if and only if it admits a translation-like action by a non-abelian free group.
\item[(3')] Geometric Gersten Conjecture: a finitely generated group is not word hyperbolic if and only if it admits a translation-like action by some Baumslag-Solitar group.
\end{enumerate}

Whyte proved (2'):

\begin{thm}[Geometric von Neumann Conjecture; Whyte, \cite{KW}] \label{INTRO THM KW}
A finitely generated group is non-amenable if and only if it admits a translation-like action by some (equivalently every) non-abelian free group.
\end{thm}

In this paper we answer (1') in the affirmative.

\begin{thm}[Geometric Burnside's Problem] \label{INTRO THM BP}
Every finitely generated infinite group admits a translation-like action by $\Z$.
\end{thm}

We actually prove something stronger than both (1') and (2'). The main distinction of the following theorem, as opposed to Theorems \ref{INTRO THM KW} and \ref{INTRO THM BP}, is the transitivity of the action.

\begin{thm} \label{TRANSITIVE}
Let $G$ be a finitely generated infinite group. Then:
\begin{enumerate}
\item[\rm (i)] $G$ has finitely many ends if and only if it admits a transitive translation-like action by $\Z$;
\item[\rm (ii)] $G$ is non-amenable if and only if it admits a transitive translation-like action by every finitely generated non-abelian free group.
\end{enumerate}
In particular, every finitely generated infinite group admits a transitive translation-like action by some (possibly cyclic) free group.
\end{thm}

In pursuit of the above theorem we prove two graph theoretic results.

\begin{thm} \label{INTRO PROP1}
Let $\Gamma$ be a connected graph whose vertices have uniformly bounded degree. Then $\Gamma$ is bilipschitz equivalent to a graph admitting a Hamiltonian path if and only if $\Gamma$ has at most two ends.
\end{thm}

\begin{thm} \label{INTRO PROP2}
Let $\Lambda_1$ and $\Lambda_2$ be two trees. If every vertex of $\Lambda_1$ and $\Lambda_2$ has degree at least three and if the vertices of $\Lambda_1$ and $\Lambda_2$ have uniformly bounded degree, then $\Lambda_1$ and $\Lambda_2$ are bilipschitz equivalent.
\end{thm}

In \cite{P}, Papasoglu proved Theorem \ref{INTRO PROP2} under the additional assumption that $\Lambda_1$ and $\Lambda_2$ are regular trees of degree at least $4$.

Theorem \ref{TRANSITIVE} above may be stated in an equivalent form without mention of translation-like actions. This alternative form, below, illustrates a peculiar property of the set of Cayley graphs associated to a group. The Cayley graphs mentioned in the following theorem are understood to be Cayley graphs coming from finite generating sets. 

\begin{thm} \label{INTRO SPAN}
If $G$ is a finitely generated infinite group then $G$ has a Cayley graph admitting a regular spanning tree. In fact, for every integer $k > 2$ the following hold:
\begin{enumerate}
\item[\rm (i)] $G$ has finitely many ends if and only if $G$ has a Cayley graph admitting a Hamiltonian path (i.e. a regular spanning tree of degree $2$);
\item[\rm (ii)] $G$ is non-amenable if and only if $G$ has a Cayley graph admitting a regular spanning tree of degree $k$.
\end{enumerate}
\end{thm}

We give an example to show that ``a Cayley graph'' cannot be replaced by ``every Cayley graph.'' Thus only certain Cayley graphs have regular spanning trees, but yet every finitely generated infinite group has at least one Cayley graph with this property.

We give a non-trivial application of Theorem \ref{INTRO SPAN} to tilings of groups, the final topic of this paper. Tilings of groups have been studied by Chou \cite{C}, Gao--Jackson--Seward \cite{GJS}, and Weiss \cite{BW} due to their usefulness in studying group actions, dynamics, equivalence relations, and general marker structures. From the algebraic and geometric viewpoints, tilings of groups are still rather mysterious as they have received little investigation and much is still unknown. In this paper we study polytiles and polytilings. These are more general notions than studied in \cite{C}, \cite{GJS}, and \cite{BW}.

\begin{defn}
For a countable group $G$, a tuple $(T_1, T_2, \ldots, T_k)$ of finite subsets of $G$, all containing the identity element, is a \emph{polytile} if there are non-empty sets $\Delta_1, \Delta_2, \ldots, \Delta_k \subseteq G$ so that $G$ is the disjoint union
$$G = \coprod_{1 \leq i \leq k, \ \delta \in \Delta_i} \delta T_i.$$
We call the tuple $P = (\Delta_1, \ldots, \Delta_k; T_1, \ldots, T_k)$ a \emph{polytiling} of $G$. The partition of $G$ expressed above is called the partition of $G$ \emph{induced} by $P$. In the case $k = 1$, polytiles are referred to as \emph{monotiles} and polytilings are referred to as \emph{monotilings}.
\end{defn}

Chou \cite{C} and Weiss \cite{BW} studied which groups $G$ have the property that for every finite $F \subseteq G$ there is a monotile $T$ with $F \subseteq T$. Weiss called groups with this property \emph{MT groups} (an acronym for ``mono-tileable''). Chou \cite{C} proved that residually finite groups, elementary amenable groups (in particular solvable groups), and groups which are free products of non-trivial groups are MT. Beyond these results, no other groups are known to be MT. Interestingly, there are also no known examples of groups which are not MT. In this paper we define and consider a related weaker question: which groups are poly-MT? Of course, care needs to be taken in defining a notion of poly-MT as any tuple $(T_1, T_2, \ldots, T_k)$ containing a singleton is a polytile. To rule out such trivialities, one wants to consider tuples $(T_1, T_2, \ldots, T_k)$ where the cardinalities of the $T_i$'s do not vary too much. The following definition seems to give the strongest restriction in this sense.

\begin{defn}
A polytile $(T_1, T_2, \ldots, T_k)$ is \emph{fair} if all of the $T_i$'s have the same number of elements. Similarly, a polytiling $(\Delta_1, \ldots, \Delta_k; T_1, \ldots, T_k)$ is \emph{fair} if the polytile $(T_1, \ldots, T_k)$ is fair.
\end{defn}

\begin{defn}
A countable group $G$ is poly-MT if for every finite $F \subseteq G$ there is a fair polytile $(T_1, \ldots, T_k)$ with $F \subseteq T_1$.
\end{defn}

We prove the following. Recall that a group is locally finite if every finite subset generates a finite subgroup.

\begin{thm} \label{INTRO POLYMT}
Every countable group is poly-MT. Furthermore, if $G$ is a countably infinite non-locally finite group then for every finite $F \subseteq G$ and all sufficiently large $n \in \N$ there is a fair polytile $(T_1, \ldots, T_k)$ with $F \subseteq T_1$ and $|T_1| = n$.
\end{thm}

In addition to studying individual tilings, as Chou and Weiss did, we also study sequences of tilings as was done by Gao--Jackson--Seward in  \cite{GJS}. The following definition generalizes definitions appearing in \cite{GJS} to the context of polytilings.

\begin{defn}
Let $G$ be a countable group, and let
$$(P_n)_{n \in \N} = (\Delta_1^n, \ldots, \Delta_{k(n)}^n; T_1^n, \ldots, T_{k(n)}^n)_{n \in \N}$$
be a sequence of polytilings. This sequence is \emph{fair} if each of the $P_n$'s is fair. The sequence is said to be \emph{coherent} if for each $n \in \N$ the partition of $G$ induced by $P_n$ is finer than the partition of $G$ induced by $P_{n+1}$. The sequence is said to be \emph{centered} if $1_G \in \Delta_1^n$ for all $n \in \N$. Finally, the sequence is \emph{cofinal} if $T_1^n \subseteq T_1^{n+1}$ for all $n$ and $G = \bigcup_{n \in \N} T_1^n$. The three adjectives ``centered, cofinal, and coherent'' are abbreviated to \emph{ccc}.
\end{defn}

We point out that if a sequence of polytilings $(P_n)_{n \in \N}$ is ccc then every two elements $g_1, g_2 \in G$ lie in the same class of the partition induced by $P_n$ for all sufficiently large $n$. Also notice that a group is MT (poly-MT) if and only if it admits a cofinal sequence of monotilings (fair polytilings). Thus having a ccc sequence of monotilings (fair polytilings) is stronger than being MT (poly-MT).

In \cite{GJS}, Gao, Jackson, and Seward studied which groups admit a ccc sequence of monotilings. They called such groups \emph{ccc groups}. They proved that the following groups are ccc groups: locally finite groups; residually finite groups; nilpotent groups; solvable groups $G$ in which $[G, G]$ is polycyclic; and groups which are free products of nontrivial groups. These are currently the only groups known to be ccc, and furthermore there are no known examples of groups which are not ccc. We consider the weaker condition of which groups are poly-ccc.

\begin{defn}
A countable group $G$ is \emph{poly-ccc} if it admits a ccc sequence of fair polytilings.
\end{defn}

Notice that just as ccc is a stronger condition than MT, poly-ccc is a stronger condition than poly-MT. We prove the following.

\begin{thm} \label{INTRO POLYCCC}
Every finitely generated group is poly-ccc. That is, every finitely generated group admits a ccc sequence of fair polytilings.
\end{thm}

Interestingly, we do not know if countable non-finitely generated groups are poly-ccc.

We point out that the proofs of these results on tilings rely heavily on the fact that every finitely generated infinite group has some Cayley graph having a regular spanning tree (Theorem \ref{INTRO SPAN}). These theorems on tilings therefore demonstrate an application of this result.

The organization of the remainder of the paper is as follows. In Section \ref{SEC PRELIM} below, we present formal definitions, notation, and some simple observations. In Section \ref{SEC HAM}, we study translation-like actions of $\Z$ which are transitive and we also study bi-infinite Hamiltonian paths on graphs. In this section we prove Theorem \ref{INTRO PROP1}. We then use these findings in Section \ref{SEC BP} to prove Theorem \ref{INTRO THM BP} (the Geometric Burnside's Problem) as well as the first clause of both Theorem \ref{TRANSITIVE} and Theorem \ref{INTRO SPAN}. Section \ref{SEC SPAN} is devoted to strengthening the Geometric von Neumann Conjecture as well as constructing Cayley graphs having regular spanning trees. In this section we prove Theorem \ref{INTRO PROP2} and prove the remaining statements of Theorem \ref{TRANSITIVE} and Theorem \ref{INTRO SPAN}. Finally, in Section \ref{SEC TILE} we study tilings of groups and prove Theorems \ref{INTRO POLYMT} and \ref{INTRO POLYCCC}.

\subsection*{Acknowledgements}
This research was supported by a National Science Foundation Graduate Research Fellowship. The author would like to thank Khalid Bou-Rabee for helpful comments regarding an earlier version of this paper.

\section{Preliminaries} \label{SEC PRELIM}

In this section we go over some basic definitions, observations, and notation. A fundamental notion which will be used frequently is that of a bilipschitz map. If $(X, d)$ and $(Y, \rho)$ are two metric spaces and $f\co X \rightarrow Y$ is a function, then $f$ is \emph{lipschitz} if there is a constant $c > 0$ such that for all $x_1, x_2 \in X$
$$\rho(f(x_1), f(x_2)) \leq c \cdot d(x_1, x_2).$$
If we additionally have that
$$\frac{1}{c} \cdot d(x_1, x_2) \leq \rho(f(x_1), f(x_2))$$
for all $x_1, x_2 \in X$, then $f$ is called \emph{bilipschitz}. The metric spaces $(X, d)$ and $(Y, \rho)$ are \emph{bilipschitz equivalent} if there is a bijective bilipschitz function $f\co X \rightarrow Y$. Finally, if $d$ and $\rho$ are two metrics on $X$, then they are called \emph{bilipschitz equivalent} if the identity map $(X, d) \rightarrow (X, \rho)$ is bilipschitz.

Notice that if $d$ and $\rho$ are bilipschitz equivalent metrics on $X$ then an action of $H$ on $X$ is translation-like with respect to the metric $d$ if and only if it is translation-like with respect to the metric $\rho$. Thus to discuss translation-like actions on $X$, it suffices to specify a family of bilipschitz equivalent metrics on $X$. This is particularly useful when $X = G$ is a finitely generated group. If $G$ is a finitely generated group, then we can specify a natural family of bilipschitz equivalent metrics on $G$. Let $S$ be any finite generating set for $G$, and let $d_S$ be the metric on $G$ given by
$$d_S(g,h) = \min\{n \in \N \: g^{-1} h \in (S \cup S^{-1})^n\}$$
(where $A^0 = \{1_G\}$ and $A^n = \{a_1 \cdot a_2 \cdots a_n \: a_i \in A\}$ for $A \subseteq G$). Notice that $d_S$ is left-invariant, meaning that $d_S(kg, kh) = d_S(g, h)$ for all $k, g, h \in G$. We call $d_S$ the \emph{left-invariant word length metric} corresponding to $S$. All of the metrics $\{d_S \: S \subseteq G \text{ is a finite generating set}\}$ are bilipschitz equivalent. To see this, one only needs to consider expressing the elements of one generating set as products of elements of the other generating set. Whenever discussing translation-like actions on a finitely generated group $G$, we will always use the family of bilipschitz equivalent metrics just specified.

Most of our arguments in this paper rely on studying the structure of various graphs. We therefore review some notation and terminology related to graphs. Let $\Gamma$ be a graph. We denote the vertex set of $\Gamma$ by $\V(\Gamma)$ and the edge set of $\Gamma$ by $\E(\Gamma)$. The graph $\Gamma$ is \emph{regular} if every vertex has the same degree, and it is \emph{even} if the degree of every vertex is finite and even. A \emph{subgraph of $\Gamma$} is a graph $\Phi$ where $\V(\Phi) \subseteq \V(\Gamma)$ and $\E(\Phi) \subseteq \E(\Gamma)$. We write $\Phi \leq \Gamma$ to denote that $\Phi$ is a subgraph of $\Gamma$. A subgraph $\Phi$ of $\Gamma$ is \emph{spanning} if $\V(\Phi) = \V(\Gamma)$. A \emph{path} $P$ in $\Gamma$ is a sequence of vertices of $\Gamma$ such that vertices which are consecutive in the sequence are joined by an edge in $\Gamma$. Formally we represent $P$ as a function from a (finite or infinite) subinterval of $\Z$ into $\V(\Gamma)$. If there is a smallest $a \in \dom(P)$ then we call $P(a)$ the \emph{initial vertex} or \emph{starting vertex} of $P$. If there is a largest $b \in \dom(P)$ then we call $P(b)$ the \emph{final vertex} or \emph{ending vertex} of $P$. If $v \in \V(\Gamma)$, we say that $P$ \emph{traverses} $v$ or \emph{visits} $v$ if there is $i \in \dom(P)$ with $P(i) = v$. If $(v_0, v_1) \in \E(\Gamma)$, we say that $P$ \emph{traverses} the edge $(v_0, v_1)$ if there is $i \in \Z$ with $i, i+1 \in \dom(P)$ and $\{P(i), P(i+1)\} = \{v_0, v_1\}$. A \emph{Hamiltonian path} is a path which visits every vertex precisely one time. These paths will be particularly important to us. A similar notion is that of an Eulerian path. An \emph{Eulerian path} is a path which traverses each edge precisely one time. We realize the vertex set of every connected graph $\Gamma$ as a metric space by using the \emph{path length metric}. Specifically, the distance between $u, v \in \V(\Gamma)$ is the infimum of the lengths of paths joining $u$ and $v$. If no metric is specified on a graph, then it is understood that we are using the path length metric. We say that two graphs are \emph{bilipschitz equivalent} if their vertex sets equipped with their respective path length metrics are bilipschitz equivalent. Finally, if $A \subseteq \V(\Gamma)$ and $k \geq 1$, we write $\rmN_k^\Gamma(A)$ to denote the union of $A$ with the set of all vertices $v \in \V(\Gamma)$ which are within distance $k$ of some vertex of $A$.

In this paper we will discuss translation-like actions not only on groups but also on connected graphs as well (using the path length metric described in the previous paragraph). Many times these graphs will arise in the form of Cayley graphs. If $G$ is a finitely generated group and $S$ is a finite generating set for $G$, then the right Cayley graph of $G$ with respect to $S$, denoted $\Cay(G; S)$, is the graph with vertex set $G$ and edge relation $\{(g, gs) \: g \in G, \ s \in S \cup S^{-1}\}$. There is a symmetric notion of a left Cayley graph where the edge set is $\{(g, sg) \: g \in G, \ s \in S \cup S^{-1}\}$. The left and right Cayley graphs are graph isomorphic via the map $g \mapsto g^{-1}$. In this paper we will always use the term \emph{Cayley graph} to mean a right Cayley graph $\Cay(G; S)$ where $S$ is a finite generating set for $G$ (we never discuss or consider Cayley graphs corresponding to infinite generating sets). If $G$ is a finitely generated group and $S \subseteq G$ is a finite generating set, then there are two corresponding metrics on $G$: $d_S$ (described earlier in this section) and the path length metric coming from $\Cay(G; S)$. It is easy to see that these two metrics are identical. Thus we can freely and without concern switch between discussing translation-like actions on $G$ and translation-like actions on a Cayley graph of $G$.

Finally, we review the notion of the number of ends of a graph and of a finitely generated group. For a connected graph $\Gamma$, the \emph{number of ends of $\Gamma$} is defined to be the supremum of the number of infinite connected components of $\Gamma - A$ as $A$ ranges over all finite subsets of $\E(\Gamma)$. If $G$ is a finitely generated group, then it is known that all Cayley graphs of $G$ have the same number of ends (see for example \cite[Section I.8, Proposition 8.29]{BH}). Thus the \emph{number of ends of $G$} is defined to be the number of ends of any Cayley graph of $G$. Finite groups have $0$ ends, and finitely generated infinite groups have either $1$, $2$, or infinitely many ends \cite[Section I.8, Theorem 8.32]{BH}.

\section{Hamiltonian paths} \label{SEC HAM}

In this section we focus on actions of $\Z$ on graphs which are both translation-like and transitive. The existence of such actions has a nice graph theoretic characterization, as the following lemma shows.

\begin{lem} \label{BILIP HAM}
A graph $\Gamma$ admits a transitive translation-like action by $\Z$ if and only if $\Gamma$ is bilipschitz equivalent to a graph admitting a bi-infinite Hamiltonian path .
\end{lem}

\begin{proof}
First suppose there is a translation-like action $*$ of $\Z$ on $\Gamma$ which is transitive. Let $d$ be the path length metric on $\Gamma$ and let $n \geq 1$ be such that $d(v, v * 1) \leq n$ for all $v \in \V(\Gamma)$ (notice that $v * 1 \neq v$, but rather $v * 0 = v$). Let $\Gamma'$ be the graph with vertex set $\V(\Gamma)$ and let there be an edge between $v_0$ and $v_1$ if and only if there is a path in $\Gamma$ of length at most $n$ joining $v_0$ to $v_1$. Let $d'$ be the path length metric on $\Gamma'$. Then for $v_0, v_1 \in \V(\Gamma)$
$$\frac{1}{n} d(v_0, v_1) \leq d'(v_0, v_1) \leq d(v_0, v_1).$$
Thus the identity map $\id\co \V(\Gamma) \rightarrow \V(\Gamma')$ is a bilipschitz bijection between $\Gamma$ and $\Gamma'$. If we fix $v \in \V(\Gamma')$ and define $P\co \Z \rightarrow \V(\Gamma')$ by $P(k) = v * k$, then $P$ is a Hamiltonian path on $\Gamma'$ (since the action of $\Z$ is free and transitive).

Now suppose that there is a graph $\Gamma'$, a Hamiltonian path $P\co \Z \rightarrow \Gamma'$, and a bilipschitz bijection $\phi\co \V(\Gamma) \rightarrow \V(\Gamma')$. First define an action of $\Z$ on $\Gamma'$ by setting
$$v * n = P(n + P^{-1}(v))$$
for $n \in \Z$ and $v \in \V(\Gamma')$. This action is free and transitive since $P$ is a Hamiltonian path. Also, it is clear that $v * n$ is at most distance $n$ from $v$. Therefore this is a transitive translation-like action of $\Z$ on $\Gamma'$. Now we define a transitive translation-like action of $\Z$ on $\Gamma$ as follows. For $n \in \Z$ and $v \in \V(\Gamma)$ set
$$v * n = \phi^{-1}(\phi(v) * n).$$
Since $\phi$ is bijective, this action is free and transitive. Also, since $\phi$ is bilipschitz and $\phi(v) * n$ is distance at most $n$ from $\phi(v)$, it follows that the distance between $v * n$ and $v$ is bounded independently of $v \in \V(\Gamma)$. Thus this is a transitive translation-like action.
\end{proof}

In the case of Cayley graphs the previous lemma takes the following more appealing form.

\begin{cor} \label{CAYLEY GRAPHS}
A finitely generated group $G$ admits a transitive translation-like action by $\Z$ if and only if $G$ has a Cayley graph admitting a bi-infinite Hamiltonian path.
\end{cor}

\begin{proof}
If $G$ admits a transitive translation-like action by $\Z$, then let $S = S^{-1}$ be a finite generating set for $G$ and let $\Cay(G; S)$ be the corresponding Cayley graph. Follow the proof of Lemma \ref{BILIP HAM} to get a graph $\Gamma'$ having a Hamiltonian path. Now notice that $\Gamma'$ is graph isomorphic to a Cayley graph $\Cay(G; T)$, where $T = S \cup S^2 \cup \cdots \cup S^n$ for some $n \in \N$. Conversely, if $\Cay(G; T)$ is a Cayley graph of $G$ which has a Hamiltonian path, then $\Cay(G; T)$ is bilipschitz equivalent to itself and so by Lemma \ref{BILIP HAM} admits a transitive translation-like action by $\Z$. This clearly is also a transitive translation-like action of $\Z$ on $G$.
\end{proof}

We now give a general sufficient condition for a graph to be bilipschitz equivalent to a graph admitting a bi-infinite Hamiltonian path.

\begin{thm}
Let $\Gamma$ be a connected infinite graph whose vertices have uniformly bounded degree. Then $\Gamma$ is bilipschitz equivalent to a graph admitting a bi-infinite Hamiltonian path if and only if $\Gamma$ has at most two ends.
\end{thm}

\begin{proof}
First suppose that $\Gamma$ is bilipschitz equivalent to $\Gamma'$ and that $\Gamma'$ has a bi-infinite Hamiltonian path. Then $\Gamma'$ contains the canonical Cayley graph of $\Z$ as a spanning subgraph. Any spanning subgraph must have at least as many ends as the original graph. So this implies that $\Gamma'$ has at most two ends. Since the number of ends of a graph is preserved by bilipschitz equivalence \cite[Section I.8, Proposition 8.29]{BH}, $\Gamma$ has at most two ends.

Now suppose that $\Gamma$ has at most two ends. Notice that $\Gamma$ cannot have $0$ ends since $\Gamma$ is infinite, connected, and locally finite. So $\Gamma$ has either $1$ end or $2$ ends. We will first construct a path on $\Gamma$ which will be very similar to a bi-infinite Eulerian path. In fact, the path we construct will visit every vertex at least once and traverse every edge at most twice. In the construction we will use a theorem of Erd\H{o}s, Gr\"{u}nwald, and Weiszfeld \cite[Section I.3, Theorem 14]{BB} which says that if $\Lambda$ is a countably infinite connected multi-graph then $\Lambda$ admits a bi-infinite Eulerian path if and only if the following conditions are satisfied:
\begin{enumerate}
\item[(i)] if the degree of a vertex is finite then it is even;
\item[(ii)] $\Lambda$ has at most $2$ ends;
\item[(iii)] if $\Phi$ is a finite, even subgraph of $\Lambda$ then $\Lambda - \E(\Phi)$ has only one infinite connected component.
\end{enumerate}
Recall that a multi-graph is a graph where loops are allowed and multiple edges between vertices are allowed.

First suppose that $\Gamma$ has $1$ end. Let $\Gamma'$ be the multi-graph obtained from $\Gamma$ by doubling each of the edges. In other words, $\V(\Gamma') = \V(\Gamma)$ and for each $u, v \in \V(\Gamma')$ $\Gamma'$ has twice as many edges joining $u$ and $v$ as $\Gamma$ does. Then $\Gamma'$ is an infinite connected one ended multi-graph and every vertex of $\Gamma'$ has even degree. Thus there is an Eulerian path $P\co \Z \rightarrow \V(\Gamma') = \V(\Gamma)$. Clearly $P$ is a path on $\Gamma$ which traverses every edge twice.

Now suppose that $\Gamma$ has two ends. We claim that there is a finite set $E \subseteq \E(\Gamma)$ such that the graph $\Gamma - E$ has precisely two connected components with both components infinite. Since $\Gamma$ has two ends, there is a finite set $F \subseteq \E(\Gamma)$ such that $\Gamma - F$ has two infinite connected components. We wish to shrink $F$ in order to connect the finite components to the infinite components, without connecting the two infinite components to one another. We do this in two steps. First, let $C$ be one of the infinite connected components of $\Gamma - F$, and let $F'$ be the set of those edges in $F$ which have an endpoint in $C$. Now $\Gamma - F'$ still has two infinite connected components, one of which is $C$. The advantage we now have is that every edge in $F'$ has an endpoint in $C$. We want to connect the finite connected components of $\Gamma - F'$ to $C$ without connecting $C$ to the other infinite connected component of $\Gamma - F'$. So for the second step let $E$ be the set of those edges in $F'$ which have both endpoints contained in infinite connected components of $\Gamma - F'$. Then $\Gamma - E$ has no finite connected components but it has two infinite connected components. This proves the claim as $E$ has the desired property.

Let $E \subseteq \E(\Gamma)$ be as in the previous paragraph. Let $C_1$ and $C_2$ be the two connected components of $\Gamma - E$. So $C_1$ and $C_2$ are infinite. Since $\Gamma$ is connected, there is $p \in C_1$ which is adjacent to $C_2$ in $\Gamma$. Let $\Lambda_1$ and $\Lambda_2$ be the induced subgraphs on $C_1$ and $C_2 \cup \{p\}$, respectively. Specifically, $\V(\Lambda_1) = C_1$, $\V(\Lambda_2) = C_2 \cup \{p\}$, and $u, v \in \V(\Lambda_i)$ are joined by an edge if and only if they are joined by an edge in $\Gamma$. For $i = 1, 2$, let $W_i : \N \rightarrow \V(\Lambda_i)$ be a path that begins at $p$ and does not self-intersect (i.e. $W_i$ is a path beginning at $p$ and going off to infinity in $\Lambda_i$). Such paths exist since $\Lambda_1$ and $\Lambda_2$ are infinite, connected, and locally finite. For $i = 1, 2$, let $\Lambda_i'$ be the multi-graph obtained from $\Lambda_i$ by doubling all edges of $\Lambda_i$ which are not traversed by $W_i$. Specifically, $\V(\Lambda_i') = \V(\Lambda_i)$ and for $u, v \in \V(\Lambda_i')$: there is no edge between $u$ and $v$ in $\Lambda_i'$ if $(u, v) \not\in \E(\Lambda_i)$; there is a single edge between $u$ and $v$ in $\Lambda_i'$ if $(u, v) \in \E(\Lambda_i)$ and $W_i$ traverses $(u, v)$; there are two edges between $u$ and $v$ in $\Lambda_i'$ if $(u, v) \in \E(\Lambda_i)$ and $W_i$ does not traverse $(u, v)$.

We claim that $p$ is the unique vertex of odd degree in $\Lambda_i'$. If $u \in \V(\Lambda_i')$ is not visited by the path $W_i$, then $u$ is connected to each of its neighbors via $2$ edges in $\Lambda_i'$ and thus $u$ has even degree in $\Lambda_i'$. If $u \neq p$ and $u$ is visited by the path $W_i$, then $W_i$ traverses precisely two distinct edges (in $\Lambda_i$) adjacent to $u$. All edges adjacent to $u$ and not traversed by $W_i$ were doubled in passing to $\Lambda_i'$. Thus $u$ has even degree in $\Lambda_i'$. Now it only remains to show that $p$ has odd degree in $\Lambda_i'$. Since only one edge adjacent to $p$ was traversed by $W_i$, all edges adjacent to $p$ but one were doubled in passing to $\Lambda_i'$. Thus $p$ has odd degree in $\Lambda_i'$, completing the proof of the claim.

We claim that the multi-graph $\Lambda_1' \cup \Lambda_2'$ satisfies the conditions (i), (ii), and (iii) of the Erd\H{o}s--Gr\"{u}nwald--Weiszfeld theorem. Notice that $p$ is the only vertex contained in both $\Lambda_1'$ and $\Lambda_2'$. Since $p$ has odd degree in each $\Lambda_i'$, it has even degree in $\Lambda_1' \cup \Lambda_2'$. So clauses (i) and (ii) are clearly satisfied. We consider (iii). In our argument we will make use of the fact that if a finite multi-graph has at most one vertex of odd degree, then it has no vertices of odd degree at all. This follows from the fact that the sums of the degrees of the vertices is always twice the number of edges. Consider a finite even subgraph $\Phi$ of $\Lambda_1' \cup \Lambda_2'$. For $i = 1, 2$ set $\Phi_i = \Phi \cap \Lambda_i'$. Since $\Phi$ is even, $\Phi_i$ has at most one vertex of odd degree, namely $p$. Since $\Phi_i$ is finite, it has no vertices of odd degree. Therefore $\Phi_i$ is a finite even subgraph of $\Lambda_i'$. Thus every vertex of $\Lambda_i' - \E(\Phi_i)$ has the same degree modulo $2$ as it has in $\Lambda_i'$. Thus $p$ is the unique vertex of odd degree in $\Lambda_i' - \E(\Phi_i)$. So any finite connected component of $\Lambda_i' - \E(\Phi_i)$ can have at most one vertex of odd degree and therefore cannot have any vertices of odd degree. So $p$ must lie in an infinite connected component of $\Lambda_i' - \E(\Phi_i)$. However, $\Gamma$ has two ends and from how we constructed $\Lambda_1'$ and $\Lambda_2'$ one can see that each $\Lambda_i'$ must be one ended. So $p$ is contained in the unique infinite connected components of $\Lambda_1' - \E(\Phi_1)$ and $\Lambda_2' - \E(\Phi_2)$.  Thus $\Lambda_1' \cup \Lambda_2' - \E(\Phi)$ has precisely one infinite connected component. By the Erd\H{o}s--Gr\"{u}nwald--Weiszfeld theorem we get an Eulerian path
$$P\co \Z \rightarrow \V(\Lambda_1') \cup \V(\Lambda_2') = \V(\Gamma).$$
Clearly $P$ is a path on $\Gamma$ which visits every vertex at least once and traverses every edge at most twice.

Regardless of whether $\Gamma$ has one end or two, we have a path $P: \Z \rightarrow \V(\Gamma)$ which visits every vertex at least once and traverses every edge at most twice. Now we will use this path $P$ to show that $\Gamma$ is bilipschitz equivalent to a graph admitting a bi-infinite Hamiltonian path. In order to do this, we will need to apply Hall's Marriage Theorem. Recall that a graph $\Lambda$ is bipartite if there is a partition $\{V_1, V_2\}$ of $\V(\Lambda)$ such that every egdge of $\Lambda$ joins a vertex of $V_1$ to a vertex of $V_2$. Hall's Marriage Theorem states that if $\Lambda$ is a locally finite bipartite graph with bipartition $\{V_1, V_2\}$ of $\V(\Lambda)$ and if $|\rmN^\Lambda_1(T) - T| \geq |T|$ for all $T \subseteq V_1$, then there is an injection $\phi\co V_1 \rightarrow V_2$ with the property that $v$ and $\phi(v)$ are joined by an edge for every $v \in V_1$ (the statement for finite graphs is Theorem 7 in Section III.3 of \cite{BB}; the theorem extends to infinite graphs by a standard compactness argument).

Let $D$ be a positive integer satisfying $\deg_\Gamma (v) \leq D$ for all $v \in \V(\Gamma)$. Set $M = D + 1$. Let $\Lambda$ be the bipartite graph with vertex set $M \Z \cup \V(\Gamma)$ and edge relation given by
$$(k \cdot M, v) \in \E(\Lambda) \Longleftrightarrow \exists 0 \leq i < M \ P(k \cdot M + i) = v.$$
Fix a finite $T \subseteq M \Z$ and set $\partial T = \rmN^\Lambda_1(T) - T \subseteq \V(\Gamma)$. We wish to show that $|\partial T| \geq |T|$. Notice that
$$\partial T = P(T) \cup P(T+1) \cup \cdots \cup P(T+M-1).$$
Let $\Phi_T \leq \Gamma$ be the graph with vertex set $\partial T$ and edge relation
$$\{(P(t+i), P(t+i+1)) \: t \in T, 0 \leq i \leq M - 2\}.$$
Since $P$ traverses each edge at most twice, $\Phi_T$ has at least $\frac{1}{2}|T|(M-1)$ edges. Since $P$ is a path in $\Gamma$, every vertex of $\Phi_T$ has degree at most $D = M-1$. Therefore
$$(M-1) |\V(\Phi_T)| \geq \sum_{v \in V(\Phi_T)} \deg_{\Phi_T} (v) = 2 \cdot |\E(\Phi_T)| \geq |T|(M-1)$$
and hence
$$|\partial T| = |\V(\Phi_T)| \geq |T|.$$
Thus the condition for Hall's Marriage Theorem is satisfied, and therefore there is a function $\phi\co M \Z \rightarrow [0, M-1]$ satisfying $(k, P(k + \phi(k))) \in \E(\Lambda)$ for all $k \in M \Z$ and $P(k_1 + \phi(k_1)) \neq P(k_2 + \phi(k_2))$ for all $k_1 \neq k_2 \in M \Z$.

Set $A = \{P(k + \phi(k)) \: k \in M \Z\} \subseteq \V(\Gamma)$. For $v \in \V(\Gamma) - A$, pick any $n_v \in \Z$ with $P(n_v) = v$. Set
$$S = \{k + \phi(k) \: k \in M \Z\} \cup \{n_v \: v \in \V(\Gamma) - A\}.$$
Clearly the restriction of $P$ to $S$ is a bijection between $S$ and $\V(\Gamma)$. Let $\psi\co \Z \rightarrow S$ be an order preserving bijection, and define $Q\co \Z \rightarrow \V(\Gamma)$ by $Q(z) = P(\psi(z))$. Then $Q$ is a bijection. By the definition of $S$, it is clear that consecutive numbers in $S$ are separated by a distance of at most $2M - 1$ and hence $\psi(k+1) - \psi(k) \leq 2M - 1$ for all $k \in \Z$. Since $P$ is a path in $\Gamma$, the distance between $Q(k+1)$ and $Q(k)$ is at most $2M - 1$. Now let $\Gamma'$ be the graph with vertex set $\V(\Gamma)$ and let there be an edge between $v_0$ and $v_1$ if and only if there is a path in $\Gamma$ of length at most $2M-1$ joining $v_0$ to $v_1$. Then $\Gamma$ and $\Gamma'$ are bilipschitz equivalent and $Q$ is a Hamiltonian path on $\Gamma'$.
\end{proof}

\begin{cor} \label{ACTION ON GRAPH}
Let $\Gamma$ be a connected infinite graph whose vertices have uniformly bounded degree. Then $\Gamma$ admits a transitive translation-like action by $\Z$ if and only if $\Gamma$ has at most two ends.
\end{cor}

\begin{proof}
This follows immediately from the previous theorem and Lemma \ref{BILIP HAM}.
\end{proof}

Among the graphs which are infinite, connected, locally finite, and have uniformly bounded degree, the above theorem completely classifies which of these graphs are bilipschitz equivalent to a graph admitting a bi-infinite Hamiltonian path. It would be interesting if nice characterizations of this property could be found among graphs which are not locally finite, or which are locally finite but whose vertices do not have uniformly bounded degree.

\begin{prob}
Find necessary and sufficient conditions for a graph to be bilipschitz equivalent to a graph admitting a bi-infinite Hamiltonian path.
\end{prob}

\section{Geometric Burnside's Problem} \label{SEC BP}

The proposition of the previous section allows us to easily resolve the Geometric Burnside's Problem.

\begin{thm}[Geometric Burnside's Problem]
Every finitely generated infinite group admits a translation-like action by $\Z$.
\end{thm}

\begin{proof}
Let $G$ be a finitely generated infinite group. If $G$ has finitely many ends, then $G$ has at most two ends \cite[Section I.8, Theorem 8.32]{BH}. Now consider any Cayley graph of $G$ and apply Corollary \ref{ACTION ON GRAPH} to get a translation-like action of $\Z$ on $G$ (notice that this action is in fact transitive). Now suppose that $G$ has infinitely many ends. Then by Stallings' Theorem \cite[Section I.8, Theorem 8.32, clause (5)]{BH}, $G$ can be expressed as an amalgamated product $A *_C B$ or HNN extension $A*_C$ with $C$ finite, $[A: C] \geq 3$, and $[B: C] \geq 2$. If $G = A *_C B$ then $ab \in G$ has infinite order whenever $a \in A - C$ and $b \in B - C$. If $G = A*_C = \langle A \cup \{t\} \ | \ \forall c \in C \ t c t^{-1} = \phi(c) \rangle$, where $\phi\co C \rightarrow A$ is an injective homomorphism, then $t$ has infinite order. Thus in either case, there is a subgroup $H \leq G$ with $H \cong \Z$. Then the natural right action of $\Z \cong H$ on $G$ is a translation-like action (see Proposition \ref{SUBGROUP ACTION}).
\end{proof}

Between the theorem above and Whyte's result in \cite{KW}, we have that Burnside's Problem and the von Neumann Conjecture both hold true in the geometric setting where subgroup containment is replaced by the existence of a translation-like action. It is interesting to wonder how much further these types of results can be taken. Two very specific questions in this direction are the following.

\begin{prob}[Whyte, \cite{KW}]
Is the Geometric Gersten Conjecture true? (See clause (3') in Section \ref{SEC INTRO} above).
\end{prob}

\begin{prob}
Does a finititely generated group have exponential growth if and only if it admits a translation-like action by a free semi-group on two generators?
\end{prob}

The algebraic version of the above problem would be: does a finitely generated group have exponential growth if and only if it contains a free sub-semi-group on two generators? Chou \cite{C} proved that every elementary amenable group of exponential growth contains a free sub-semi-group on two generators. However, the algebraic version is false in general. Grigorchuk has pointed out that the wreath product of the Grigorchuk group with the cyclic group of order $2$ is a torsion group of exponential growth. Since it is torsion, it cannot contain a free sub-semi-group. The geometric (translation-like action) version may hold true. By Chou's result, it holds for elementary amenable groups, and by the Geometric von Neumann Conjecture it holds for non-amenble groups. 

A more general, open-ended question is the following.

\begin{prob}
What other algebraic/geometric properties can be reformulated or characterized in terms of translation-like actions?
\end{prob}

Before ending this section we present two corollaries treating Hamiltonian paths on Cayley graphs and translation-like actions of $\Z$ which are transitive.

\begin{cor}\label{TRANS BURNS}
A finitely generated infinite group has finitely many ends if and only if it admits a translation-like action by $\Z$ which is transitive.
\end{cor}

\begin{proof}
If $G$ has finitely many ends then the proof of the previous theorem shows that there is a translation-like action of $\Z$ on $G$ which is transitive. Now suppose that there is a translation-like action of $\Z$ on $G$ which is transitive. If $\Cay(G)$ is any Cayley graph of $G$ then the action of $\Z$ on $G$ gives rise to a transitive translation-like action of $\Z$ on $\Cay(G)$. By Corollary \ref{ACTION ON GRAPH}, $\Cay(G)$ has at most two ends and thus $G$ has at most two ends.
\end{proof}

Notice that the following corollary no longer requires the group to be infinite.

\begin{cor}\label{CAY HAM}
A finitely generated group has finitely many ends if and only if it has a Cayley graph admitting a Hamiltonian path.
\end{cor}

\begin{proof}
For infinite groups, this follows from the above corollary and Corollary \ref{CAYLEY GRAPHS}. Finally, every finite group has $0$ ends (in particular, has finitely many ends), and the complete graph on $|G|$ vertices is a Cayley graph of $G$ which clearly admits a Hamiltonian path.
\end{proof}

Relating to the above corollary, we point out that Igor Pak and Rado\v{s} Radoi\v{c}i\'c proved in \cite{PR} that every finite group $G$ has a generating set $S$ with $|S| \leq \log_2 |G|$ such that the Cayley graph $\Cay(G; S)$ admits a Hamiltonian path. This ties in with the well known Lov\'asz Conjecture: every finite, connected, vertex-transitive graph admits a Hamiltonian path. Clearly this conjecture applies to finite Cayley graphs. Despite the Pak--Radoi\v{c}i\'c result, the special case of the Lov\'asz Conjecture for Cayley graphs is far from settled.

\begin{prob} [Special case of the Lov\'asz Conjecture]
Does every Cayley graph of every finite group admit a Hamiltonian path?
\end{prob}

\begin{prob}
Does every Cayley graph of every finitely generated group with finitely many ends admit a Hamiltonian path?
\end{prob}

\section{Regular spanning trees of Cayley graphs} \label{SEC SPAN}

In this section we will prove that every finitely generated infinite group has a Cayley graph containing a regular spanning tree. To do this we will strengthen Whyte's result on the Geometric von Neumann Conjecture.

We first present two simple properties of translation-like actions. These properties were mentioned in \cite{KW} without proof. For convenience to the reader we include a proof here.

\begin{lem}[Whyte, \cite{KW}] \label{TRANS ORB LIP}
Let $H$ be a finitely generated group, let $(X, d)$ be a metric space, and suppose $H$ acts on $X$ and the action is translation-like. Then there is a constant $C$ such that for all $x \in X$ the map
$$h \in H \mapsto x * h \in X$$
is a $C$--lipschitz injection.
\end{lem}

\begin{proof}
Let $U$ be a finite generating set for $H$ and let $\rho$ be the corresponding left-invariant word-length metric. For each $u \in U \cup U^{-1}$ let $c_u \geq 1$ be such that $d(x, x * u) \leq c_u$ for all $x \in X$. Set $C = \max\{c_u \: u \in U \cup U^{-1}\}$ and fix $x \in X$. For $h, k \in H$ we have
$$d(x * k, x * h) = d( (x * h) * h^{-1} k, x * h) \leq C \cdot \rho(h^{-1} k, 1_H) = C \cdot \rho(k, h).$$
So the map is $C$--lipschitz as claimed. Also, the map is injective as translation-like actions are free.
\end{proof}

\begin{cor}[Whyte, \cite{KW}] \label{GRAPH PARTITION}
Let $H$ and $G$ be finitely generated groups, let $\Cay(H)$ be a Cayley graph of $H$, and suppose that $H$ acts on $G$ and the action is translation-like. Then there is a Cayley graph $\Cay(G)$ of $G$ so that for all $g \in G$ and $h_0, h_1 \in H$ with $(h_0, h_1) \in \E(\Cay(H))$, $(g * h_0, g * h_1) \in \E(\Cay(G))$. In particular, there is a spanning subgraph $\Phi$ of $\Cay(G)$ such that the connected components of $\Phi$ are precisely the orbits of the $H$ action and every connected component of $\Phi$ is graph isomorphic to $\Cay(H)$.
\end{cor}

\begin{proof}
Let $V = V^{-1}$ be a finite generating set for $G$, and let $d$ be the corresponding left-invariant word length metric. Let $\Cay(H)$ be a Cayley graph of $H$, say $\Cay(H) = \Cay(H; U)$ where $U = U^{-1}$ is a finite generating set for $H$. Let $\rho$ be the left-invariant word length metric on $H$ corresponding to $U$. Let $C$ be as in the previous lemma. By picking a larger $C$ if necessary, we may assume $C \in \N$. Set $W = V \cup V^2 \cup \cdots \cup V^C$. Then $W$ is a generating set for $G$ and for all $g \in G$ and $u \in U \subseteq H$ we have $g *u \in g W$. Consider the Cayley graph $\Cay(G; W)$. Fix $g \in G$ and $h_0, h_1 \in H$ with $(h_0, h_1) \in \E(\Cay(H;U))$. Then there is $u \in U$ with $h_1 = h_0 u$. Therefore
$$g * h_1 = (g * h_0) * u \in (g * h_0) W,$$
so $(g * h_1, g * h_0) \in \E(\Cay(G; W))$. Thus the first claim holds. Now let $\Phi$ be the graph with vertex set $G$ and edge set $\{(g, g * u) \: g \in G, \ u \in U\}$. Then $\Phi$ is a spanning subgraph of $\Cay(G; W)$, the connected components of $\Phi$ are the orbits of the $H$ action, and every connected component of $\Phi$ is graph isomorphic to $\Cay(H; U)$ (since the action of $H$ is free).
\end{proof}

Let $\Lambda$ be a tree with a distinguished root vertex $\lambda$. Define a partial ordering, denoted $\prec$, on $\V(\Lambda)$ by declaring $u \prec v$ if and only if the unique shortest path from $\lambda$ to $v$ traverses $u$.

\begin{defn}
We call a finite set $P \subseteq \V(\Lambda)$ a \emph{perimeter} if $P \neq \{\lambda\}$ and the following conditions are met:
\begin{enumerate}
\item[(i)] there is a constant $R$ such that if $d(u, \lambda) \geq R$ then there is $p \in P$ with $p \prec u$;
\item[(ii)] if $p, p' \in P$ and $p \prec p'$, then $p = p'$.
\end{enumerate}
The smallest such $R$ satisfying (i) is the \emph{radius of $P$}.
\end{defn}

Perimeters will soon be used to help us construct quasi-isometries between trees. The idea is that if $\Lambda$ and $T$ are trees with distinguished roots, $F\co \V(T) \rightarrow \V(\Lambda)$ is a partially defined function, and $F(t)$ is defined, then we want to extend $F$ so that the vertices of $T$ adjacent to $t$ are mapped, not necessarily injectively, to a perimeter of $F(t)$. The following is a technical lemma needed for the proposition following it where this construction is carried out. From the viewpoint of the construction just discussed, in the lemma below $r$ essentially represents the number of vertices of $T$ adjacent to $t$, and $d(p)$ tells us how many vertices of $T$ adjacent to $t$ should be mapped to $p$.

\begin{lem} \label{LEM PERMTR}
Let $\Lambda$ be a tree with distinguished root vertex $\lambda$, let $\prec$ be defined as above, and let $r \geq 3$. If $2 \leq \deg(\lambda) \leq r - 1$ and all the other vertices of $\Lambda$ have degree at least $3$, then there is a perimeter $P$ of radius at most $r$ and a function $d\co P \rightarrow \N^+$ satisfying $d(p) \leq \deg(p) - 1$ for all $p \in P$ and
$$\sum\limits_{p \in P} d(p) = r.$$
\end{lem}

\begin{proof}
The proof is by induction on $r$. If $r = 3$, then $\deg(\lambda) = 2$. Let $P$ be the set of vertices adjacent to $\lambda$, say $P = \{a, b\}$. Define $d(a) = 2$ and $d(b) = 1$. Then $P$ and $d$ have the desired properties. Now suppose that $r \geq 4$ and that the claim is true for all $r' < r$. We have $\deg(\lambda) \leq r-1$, so by the Euclidean algorithm there are $q, s \in \N$ with $s < \deg(\lambda)$ and
$$r = q \cdot \deg(\lambda) + s = q( \deg(\lambda) - s) + (q + 1)s.$$
Clearly $q \leq \frac{r}{\deg(\lambda)} \leq \frac{r}{2}$. Since $r \geq 4$, we have $q, q+1 < r$.

Divide the vertices adjacent to $\lambda$ into two disjoint sets $A$ and $B$ with $|A| = \deg(\lambda) - s > 0$ and $|B| = s$ (note that we can have $s = 0$ in which case $B = \varnothing$). Fix a vertex $a \in A$. If $q \leq \deg(a) - 1$, then set $P_a = \{a\}$ and $d_a(a) = q$. If $q > \deg(a) - 1$, then let $\Lambda_a$ be the graph with vertex set $\{u \in \V(\Lambda) \: a \prec u\}$ and define there to be an edge between $u$ and $v$ if and only if they are joined by an edge in $\Lambda$. Then $\Lambda_a$ is a tree, and we declare $a$ to be its root vertex. Clearly $\deg_{\Lambda_a}(a) = \deg_{\Lambda}(a) - 1$ and $\deg_{\Lambda_a}(u) = \deg_\Lambda(u)$ for all $a \neq u \in \V(\Lambda_a)$. By the inductive hypothesis, there is a perimeter $P_a$ of $\Lambda_a$ of radius at most $q$ and a function $d_a\co P_a \rightarrow \N^+$ with $d_a(p) \leq \deg_{\Lambda_a}(p) - 1 = \deg_{\Lambda}(p) - 1$ for all $p \in P_a$ and
$$\sum\limits_{p \in P_a} d_a(p) = q.$$
Carry out the same construction for every element of $B$, but with $q$ replaced by $q+1$. Set
$$P = \bigcup_{c \in A \cup B} P_c \text{ and } d = \bigcup_{c \in A \cup B} d_c.$$
Then $P$ is a perimeter for $\Lambda$, $d(p) \leq \deg_{\Lambda}(p) - 1$ for all $p \in P$, and
$$\sum\limits_{p \in P} d(p) = \sum\limits_{a \in A} \sum\limits_{p \in P_a} d_a(p) + \sum\limits_{b \in B} \sum\limits_{p \in P_b} d_b(p) = q |A| + (q+1) |B| = r.$$
Also, the radius of $P$ is at most
$$q + 2 \leq \frac{r}{2} + 2 \leq r.$$
\end{proof}

\begin{thm}\label{TWO TREES}
Let $\Lambda_1$ and $\Lambda_2$ be two trees. If every vertex of $\Lambda_1$ and $\Lambda_2$ has degree at least three and if the vertices of $\Lambda_1$ and $\Lambda_2$ have uniformly bounded degree, then $\Lambda_1$ and $\Lambda_2$ are bilipschitz equivalent.
\end{thm}

\begin{proof}
Let $r > 2$ be such that every vertex of $\Lambda_1$ and $\Lambda_2$ has degree at most $r$. Let $T$ be the regular tree of degree $r+1$. It suffices to show that $\Lambda_1$ and $\Lambda_2$ are both bilipschitz equivalent to $T$. So dropping subscripts, let $\Lambda$ be a tree such that every vertex has degree at least $3$ and at most $r$. We will show that $\Lambda$ and $T$ are bilipschitz equivalent.

We will first define a quasi-isometry $F\co \V(T) \rightarrow \V(\Lambda)$. Recall that a function $f\co (X, d) \rightarrow (Y, \rho)$ between metric spaces $(X, d)$ and $(Y, \rho)$ is called a \emph{quasi-isometry} if there are constants $K$ and $C$ such that
$$\frac{1}{K} d(x_1, x_2) - C \leq \rho(f(x_1), f(x_2)) \leq K d(x_1, x_2) + C$$
and such that for every $y \in Y$ there is $x \in X$ with $\rho(y, f(x)) < C$. Also recall that $X$ and $Y$ are \emph{quasi-isometric} if there is a quasi-isometry $f\co X \rightarrow Y$. We will view $\V(T)$ and $\V(\Lambda)$ as metric spaces by using the path length metrics described in Section \ref{SEC PRELIM}. Let $d$ and $\rho$ denote the path length metrics on $\V(T)$ and $\V(\Lambda)$, respectively. Fix root vertices $t_0 \in \V(T)$ and $\lambda_0 \in \V(\Lambda)$. Define a partial ordering, denoted $\prec$, on $\V(T)$ by declaring $u \prec v$ if and only if the unique shortest path from $t_0$ to $v$ traverses $u$. Similarly define a partial ordering, also denoted $\prec$, on $\V(\Lambda)$. Define
$$S(T; n) = \{u \in \V(T) \: d(t_0, u) = n\} \text{ and } B(T; n) = \{u \in \V(T) \: d(t_0, u) \leq n\}.$$
Clearly $\V(T)$ is the union of the $S(T; n)$'s. We will inductively define a quasi-isometry $F\co \V(T) \rightarrow \V(\Lambda)$ so that the following properties hold:
\begin{enumerate}
\item[(i)] $F(t_0) = \lambda_0$;
\item[(ii)] if $F(u) = F(v)$, then $d(u, v) \leq 2$ and $d(t_0, u) = d(t_0, v)$;
\item[(iii)] if $d(u, v) = 1$ then $\rho(F(u), F(v)) \leq r+1$;
\item[(iv)] if $u \prec v$ then $F(u) \prec F(v)$;
\item[(v)] if $F(u) \prec F(v)$ then there is $w \in F^{-1}(u)$ with $w \prec v$;
\item[(vi)] for $y \in \V(\Lambda)$, $|F^{-1}(y)| \leq \deg_\Lambda(y) - 1$;
\item[(vii)] if $u, v \in S(T; n)$ and $F(u) \prec F(v)$, then $F(u) = F(v)$;
\item[(viii)] for every $n \geq 1$ $F(S(T; n))$ is a perimeter.
\end{enumerate}

Set $F(t_0) = \lambda_0$ so that (i) is satisfied. By Lemma \ref{LEM PERMTR}, there is a perimeter $P$ of $\Lambda$ of radius at most $r+1$ and a function $d\co P \rightarrow \N^+$ with $d(p) \leq \deg_\Lambda(p) - 1$ for all $p \in P$ and $\sum_{p \in P} d(p) = r+1$. Since $|S(T;1)| = r+1$, there is a surjection $F\co S(T;1) \rightarrow P$ such that $|F^{-1}(p)| = d(p)$. The definition of $F$ on $B(T; 1)$ satisfies clauses (i) through (viii).

Now suppose that $F$ has been defined on $B(T; n-1)$ and does not violate any of the clauses (i) through (viii). We will define $F$ on $S(T; n)$. Notice that (v) and (ii) tell us that if $F(u) \prec F(v)$ then $d(t_0, u) \leq d(t_0, v)$. Therefore by (vii) we have that if $u \in S(T; n-1)$ then there is nothing strictly $\prec$--larger than $F(u)$ in $F(B(T; n-1))$. This together with (viii) imply that the $\prec$--maximal elements of $F(B(T; n-1))$ are precisely the $F$ images of elements of $S(T; n-1)$. Let $y_1, y_2, \ldots, y_m$ be the $\prec$--maximal elements of $F(B(T;n-1))$, and set $U_i = F^{-1}(y_i) \cap S(T; n-1) \neq \varnothing$. Then the $U_i$'s form a partition of $S(T; n-1)$. To define $F$ on $S(T;n)$, it suffices to fix an $i$ and define $F$ on the elements of $S(T;n)$ which are adjacent to an element of $U_i$. This definition is provided in the following paragraph. To simplify notation, in the paragraph below we work with a single $y_i$ and $U_i$ but omit the subscript $i$.

By (vi), $|U| \leq \deg_\Lambda(y) - 1$. Fix any function $h\co U \rightarrow \N^+$ satisfying $\sum_{v \in U} h(v) = \deg_\Lambda(y) - 1$. Partition the set of vertices which are adjacent to $y$ and $\prec$--larger than $y$ into disjoint sets $\{A_v \: v \in U\}$ such that $|A_v| = h(v)$. For each $v \in U$ we will define a graph $\Lambda_v$. We first describe the vertex set. If $|A_v| \geq 2$, then we let $\Lambda_v$ be the graph with vertex set the union of $y$ with $A_v$ and with all of the vertices of $\Lambda$ which are $\prec$--larger than some element of $A_v$. If $|A_v| = 1$ then we let $\Lambda_v$ have the same vertices as described in the previous sentence, except that we exclude $y$ from the vertex set. In either case, we let there be an edge between $w_1, w_2 \in \V(\Lambda_v)$ if and only if they are joined by an edge in $\Lambda$. Notice that $\Lambda_v$ is a subgraph of $\Lambda$, $\deg_{\Lambda_v}(y) = |A_v| \leq \deg_\Lambda(y) - 1$ (if $y \in \V(\Lambda_v)$), $\deg_{\Lambda_v}(a) = \deg_\Lambda(a) - 1$ (if $A_v = \{a\}$), and $\deg_{\Lambda_v}(w) = \deg_\Lambda(w)$ for all other vertices $w \in \V(\Lambda_v)$. We declare the root vertex of $\Lambda_v$ to be $y$ if $y \in \V(\Lambda_v)$ and otherwise to be $a$ if $A_v = \{a\}$. By Lemma \ref{LEM PERMTR}, there is a perimeter $P_v$ of $\Lambda_v$ of radius at most $r$ and a function $d_v\co P_v \rightarrow \N^+$ with $d_v(p) \leq \deg_{\Lambda_v}(p) - 1 = \deg_{\Lambda}(p) - 1$ for all $p \in P_v$ and
$$\sum\limits_{p \in P_v} d_v(p) = r.$$
Thus, if $S(v)$ denotes the set of vertices of $T$ adjacent to $v$ and $\prec$--larger than $v$, then $|S(v)| = r$ and there is a surjection $F\co S(v) \rightarrow P_v$ with $|F^{-1}(p)| = d_v(p)$ for all $p \in P_v$. This defines $F$ on $S(v)$ for each $v \in U$. Thus we have defined $F$ on all elements of $S(T; n)$ which are adjacent to some element of $U$. By varying the set $U$, we get a definition of $F$ on $S(T;n)$ and thus on $B(T;n)$. It is easy to check that $F$ still satisfies clauses (i) through (viii). By induction, this completes the definition of $F$.

We now verify, relying only on clauses (i) through (viii), that $F$ is a quasi-isometry. Clearly clause (iii) implies that $\rho(F(u), F(v)) \leq (r+1) d(u, v)$ for all $u, v \in \V(T)$. We now check by cases that for all $u, v \in \V(T)$ we have $d(u, v) - r - 2 \leq \rho(F(u), F(v))$. Fix $u, v \in \V(T)$ and let $w \in \V(T)$ be the $\prec$--largest vertex with $w \prec u, v$. Notice that clauses (ii) and (iv) imply that if $p \prec q$ then $d(p, q) \leq \rho(F(p), F(q))$.

\underline{Case 1:} $d(w, u) \leq 1$ or $d(w, v) \leq 1$. Without loss of generality, suppose that $d(u, w) \leq 1$. Then $\rho(F(u), F(w)) \leq r+1$  by (iii) and $d(w, v) \leq \rho(F(w), F(v))$ since $w \prec v$. So
$$d(u, v) - r - 2 = d(u, w) + d(w, v) - r - 2$$
$$\leq \rho(F(w), F(v)) - \rho(F(w), F(u)) \leq \rho(F(u), F(v)).$$

\underline{Case 2:} $d(w, u), d(w, v) \geq 2$. Let $w_u$ and $w_v$ be the unique vertices of $T$ with $d(w, w_u) = d(w, w_v) = 2$, $w \prec w_u \prec u$, and $w \prec w_v \prec v$. Notice that $d(w_u, w_v) = 4$. By (v) and (ii) we have that $F(w_u)$ and $F(w_v)$ are not $\prec$--comparable. Thus $\rho(F(w_u), F(w_v)) \geq 2$, and since $w_u \prec u$ and $w_v \prec v$ we have that $d(w_u, u) \leq \rho(F(w_u), F(u))$ and $d(w_v, v) \leq \rho(F(w_v), F(v))$. Thus
$$d(u, v) - 2 = d(u, w_u) + d(v, w_v) + d(w_u, w_v) - 2$$
$$\leq \rho(F(u), F(w_u)) + \rho(F(v), F(w_v)) + \rho(F(w_u), F(w_v)) = \rho(F(u), F(v)).$$

To show that $F$ is a quasi-isometry, all that remains is to show that for all $y \in \V(\Lambda)$ there is $u \in \V(T)$ with $\rho(y, F(u)) \leq r+1$. So fix $y \in \V(\Lambda)$ and towards a contradiction suppose that every image point of $F$ is at least a distance of $r+2$ from $y$. First suppose that there is $u \in \V(T)$ with $y \prec F(u)$. Let $P$ be a path in $T$ from $t_0$ to $u$. Let $n$ be least such that $y \not\prec F(P(n))$ and $y \prec F(P(n+1))$. Since $F(P(n))$ and $F(P(n+1))$ are distance at least $r+2$ from $y$, we have $\rho(F(P(n)), F(P(n+1))) \geq 2r+4$, contradicting (iii). So there are no $u \in \V(T)$ with $y \prec F(u)$. Since each $F(S(T; n))$ is a perimeter, it follows that for each $n \in \N$ there is $u_n \in S(T; n)$ with $F(u_n) \prec y$. However, there are only finitely many $w \in \V(\Lambda)$ with $w \prec y$, so this violates clause (vi). We conclude that every vertex of $\Lambda$ is within a distance of $r+1$ of some image point of $F$. Thus $F$ is a quasi-isometry, and so $\Lambda$ is quasi-isometric to $T$.

Recall that a graph $\Gamma$ is \emph{non-amenable} if there is $C > 1$ and $k \in \N$ so that $|\rmN_k^\Gamma(S)| \geq C \cdot |S|$ for all finite $S \subseteq \V(\Gamma)$, where $\rmN_k^\Gamma(S)$ is the union of $S$ with all vertices of $\Gamma$ within distance $k$ of $S$. When $\Gamma$ is the Cayley graph of a group, this notion coincides with the standard notion of amenability of groups. It is known that if $\Gamma$ is a tree and every vertex of $\Gamma$ has degree at least $3$, then $\Gamma$ is non-amenable (one can use $C = 2$ and $k = 1$). By Complement 46 of Section IV.B on page 104 of \cite{PH}, any two non-amenable quasi-isometric graphs are bilipschitz equivalent. Since $\Lambda$ and $T$ are non-amenable and quasi-isometric, they are bilipschitz equivalent.
\end{proof}

In the following theorem we strengthen Whyte's result on the Geometric von Neumann Conjecture. In the remainder of this section, we let $\F_k$ denote the non-abelian free group of rank $k$ for $k > 1$.

\begin{thm}[Geometric von Neumann Conjecture]\label{VN CONJ}
Let $G$ be a finitely generated group. The following are equivalent for every integer $k \geq 2$:
\begin{enumerate}
\item[\rm (i)] $G$ is non-amenable;
\item[\rm (ii)] $G$ admits a translation-like action by $\F_k$;
\item[\rm (iii)] $G$ admits a transitive translation-like action by $\F_k$.
\end{enumerate}
\end{thm}

\begin{proof}
In \cite{KW}, Whyte proved the equivalence of (i) and (ii). Clearly (iii) implies (ii), so we only must show that (ii) implies (iii). Let $\Cay(\F_k)$ be the canonical Cayley graph of $\F_k$. So $\Cay(\F_k)$ is a regular tree of degree $2k$. By assumption, there is a translation-like action $*$ of $\F_k$ on $G$. By Corollary \ref{GRAPH PARTITION}, there is a Cayley graph $\Cay(G)$ of $G$ having a spanning subgraph $\Phi$ in which the connected components of $\Phi$ are the orbits of the $\F_k$ action and every connected component of $\Phi$ is graph isomorphic to $\Cay(\F_k)$, the regular tree of degree $2k$.

Let $\Psi$ be the quotient graph of $\Cay(G)$ obtained by identifying points which lie in a common $\F_k$ orbit. Specifically, the vertex set of $\Psi$ is $\{g * \F_k \: g \in G\}$ and there is an edge between $g * \F_k$ and $h * \F_k$ if and only if there are $t_0, t_1 \in \F_k$ such that $g * t_0$ and $h * t_1$ are joined by an edge in $\Cay(G)$. Since $\Cay(G)$ is connected, so is $\Psi$. We construct a spanning tree $\Psi'$ of $\Psi$ as follows. Set $\psi_0 = 1_G * \F_k \in \V(\Psi)$, and let $d$ be the path length metric on $\Psi$. For each $\psi_0 \neq v \in \V(\Psi)$, pick any $f(v) \in \V(\Psi)$ with the property that $d(v, f(v)) = 1$ and $d(\psi_0, f(v)) = d(\psi_0, v) - 1$. We set $\V(\Psi') = \V(\Psi)$ and
$$\E(\Psi') = \{(v, f(v)) \: \psi_0 \neq v \in \V(\Psi)\} \cup \{(f(v), v) \: \psi_0 \neq v \in \V(\Psi)\}.$$
A simple inductive argument on the magnitude of $d(\psi_0, v)$ shows that $\Psi'$ is connected. Also, if $u$ and $v$ are joined by an edge in $\Psi'$ then $d(\psi_0, u) \neq d(\psi_0, v)$ and $d(\psi_0, u) < d(\psi_0, v)$ implies $u = f(v)$. Thus $\Psi'$ does not contain any cycles as cycles cannot contain a vertex of maximal distance to $\psi_0$. So $\Psi'$ is a spanning tree of $\Psi$.

Recall that $\Phi$ is a subgraph of $\Cay(G)$ in which the connected components are regular trees of degree $2k$. Also notice that the vertices of $\Psi$ and $\Psi'$ are precisely the connected components of $\Phi$ (since both are just the orbits of the $\F_k$ action). We will use $\Psi'$ to enlarge $\Phi$ to a spanning tree $\Phi'$ of $\Cay(G)$. Let $e \in \E(\Psi')$ be an edge in $\Psi'$, say between $g_1 * \F_k$ and $g_2 * \F_k$. Then by definition, there are $t_1, t_2 \in \F_k$ such that $g_1 * t_1$ and $g_2 * t_2$ are joined by an edge in $\Cay(G)$. Let $h(e)$ be any such edge $(g_1 * t_1, g_2 * t_2) \in \E(\Cay(G))$. So $h\co \E(\Psi') \rightarrow \E(\Cay(G))$. Now let $\Phi'$ be the graph with vertex set $G$ and edge set $h(\E(\Psi')) \cup \E(\Phi)$. Since the connected components of $\Phi$ are trees, the vertices of $\Psi'$ are the connected components of $\Phi$, and $\Psi'$ is a tree, it is easy to see that $\Phi'$ is also a tree. So $\Phi'$ is a spanning tree of $\Cay(G)$. Each vertex of $\Phi'$ has degree at least $2k$, and the degrees of the vertices of $\Phi'$ are uniformly bounded since $\Phi \leq \Phi' \leq \Cay(G)$. By Theorem \ref{TWO TREES}, $\Phi'$ (with its path length metric) is bilipschitz equivalent to the $2k$ regular tree $\Cay(\F_k)$. Let $P\co \F_k \rightarrow \V(\Phi')$ be this bilipschitz equivalence. Since $\Phi'$ is a subgraph of $\Cay(G)$, it follows that the map $P\co \F_k \rightarrow G$ is lipschitz when $G = \V(\Phi')$ is equipped with the path length metric of $\Cay(G)$. Now we define a transitive translation-like action $\circ$ of $\F_k$ on $G$ as follows. For $t \in \F_k$ and $g \in G$, we define
$$g \circ t = P(P^{-1}(g) \cdot t),$$
where $\cdot$ denotes the group operation of $\F_k$. This is clearly a transitive translation-like action since $P\co \F_k \rightarrow G$ is a lipschitz bijection (where $G$ has the path length metric coming from $\Cay(G)$).
\end{proof}

We point out that in the previous proof, instead of using Whyte's result we could have also used a similar result by Benjamini--Schramm \cite[Remark 1.9]{BS}.

Combining the previous theorem with Corollary \ref{TRANS BURNS} gives the following.

\begin{cor} \label{COR TRANS ACT FREE}
Every finitely generated infinite group admits a transitive translation-like action by some (possibly cyclic) free group.
\end{cor}

\begin{proof}
By Corollary \ref{TRANS BURNS}, $G$ admits a transitive translation-like action by $\Z$ (which is a cyclic free group) if $G$ has finitely many ends. So now suppose that $G$ has infinitely many ends. We will show that $G$ contains a non-abelian free group of rank $2$ and is thus non-amenable. By Stallings' Theorem \cite[Section I.8, Theorem 8.32, clause (5)]{BH}, $G$ can be expressed as an amalgamated product $A *_C B$ or HNN extension $A*_C$ with $C$ finite, $[A: C] \geq 3$, and $[B: C] \geq 2$. We now show by cases that $G$ contains a non-abelian free group of rank $2$.

\underline{Case 1:} $C$ is trivial and $G$ is of the form $A *_C B = A * B$. By Proposition 4 on page 6 of \cite{S}, the kernel of the homomorphism $A*B \rightarrow A \times B$ is a free group with free basis $\{a^{-1} b^{-1} a b \: 1_A \neq a \in A, \ 1_B \neq b \in B\}$. So we only need to show that this free basis contains more than one element. Since $|A| \geq 3$, there are non-identity $a_1 \neq a_2 \in A$. Fix any non-identity $b \in B$. Since $G = A * B$, it is clear that $a_1^{-1} b^{-1} a_1 b$ and $a_2^{-1} b^{-1} a_2 b$ are distinct. Thus the free basis of the kernel contains at least two elements and therefore $G$ contains a non-abelian free group of rank $2$.

\underline{Case 2:} $C$ is non-trivial and $G$ is of the form $A *_C B$. We will use Lemma 6.4 of Section III.$\Gamma$.6 on page 498 of \cite{BH}. This lemma states, in particular, that if $a_1 \in A$, $a_2, a_3, \ldots, a_n \in A - C$, $b_1, \ldots, b_{n-1} \in B - C$, and $b_n \in B$, then $a_1 b_1 a_2 b_2 \cdots a_n b_n$ is not the identity element in $G = A *_C B$. Pick any $b_1, b_2 \in B - C$, and pick any $a_1 \in A - C$. Since $[A : C] \geq 3$, we can also pick $a_2 \in A - C$ with $a_2 a_1 \not\in C$. If $b_2 b_1 \in C$, then pick $a_3 \in A - C$ with $a_3 b_2 b_1 a_1 \not\in C$. If $b_2 b_1 \not\in C$, then pick any $a_3 \in A - C$. Set $u = a_1 b_1 a_2$ and $v = b_1 a_1 b_2 a_3 b_2$. Then $u$ and $v$ generate a non-abelian free group of rank $2$.

\underline{Case 3:} $G$ is of the form $A*_C = \langle A, t | t^{-1} c t = \phi(c)\}$ where $[A : C] \geq 3$ and $\phi\co C \rightarrow A$ is an injective homomorphism. Let $C' = \phi(C)$. We again use Lemma 6.4 of Section III.$\Gamma$.6 on page 498 of \cite{BH}. This lemma states, in particular, that if $a_1, a_2, \ldots, a_{n-1} \in A - (C \cup C')$ and $m_1, m_2, \ldots, m_n \in \Z - \{0\}$ then the product $t^{m_1} a_1 t^{m_2} a_2 \cdots t^{m_{n-1}} a_{n-1} t^{m_n}$ is not the identity element in $G = A*_C$. Since $C$ is finite, $C'$ cannot properly contain $C$ and therefore $C'$ cannot contain any coset $a C$ with $a \not\in C$. Thus there is $a \in A - C \cup C'$. Set $u = t a t$ and $v = t^2 a t^2$. Then $u$ and $v$ generate a non-abelian free group of rank $2$.

Thus $G$ contains a non-abelian free group of rank $2$. The group $G$ is therefore non-amenable (this is a well known consequence of containing a non-abelian free group but also follows directly from the previous theorem), and so by the previous theorem there is a transitive translation-like action of the non-abelian free group of rank $2$ on $G$.
\end{proof}

The previous Corollary has an interesting reinterpretation in terms of Cayley graphs.

\begin{cor} \label{COR SPAN TREE}
If $G$ is a finitely generated infinite group then $G$ has a Cayley graph admitting a regular spanning tree. In fact, for every integer $k > 2$ the following hold:
\begin{enumerate}
\item[\rm (i)] $G$ has finitely many ends if and only if $G$ has a Cayley graph admitting a Hamiltonian path (i.e. a regular spanning tree of degree $2$);
\item[\rm (ii)] $G$ is non-amenable if and only if $G$ has a Cayley graph admitting a regular spanning tree of degree $k$.
\end{enumerate}
\end{cor}

\begin{proof}
For the first statement, apply the previous corollary and Corollary \ref{GRAPH PARTITION} and notice that the graph $\Phi$ referred to in Corollary \ref{GRAPH PARTITION} has only one connected component since the action is transitive. Clause (i) is Corollary \ref{CAY HAM}.

Now fix $k > 2$. If $G$ is non-amenable, then by clause (iii) of Theorem \ref{VN CONJ} and Corollary \ref{GRAPH PARTITION} $G$ has a Cayley graph $\Cay(G)$ admitting a regular spanning tree $\Phi$ of degree $2k$. Say $\Cay(G) = \Cay(G; V)$ where $V$ is a finite generating set for $G$, and let $d$ be the corresponding left-invariant word length metric. Let $\Lambda$ be the regular tree of degree $k$. By Theorem \ref{TWO TREES}, $\Phi$ and $\Lambda$ are bilipschitz equivalent, say via $f\co \V(\Lambda) \rightarrow \V(\Phi)$. Since $\Phi$ is a subgraph of $\Cay(G)$, it follows that there is a constant $c \in \N$ so that $d(f(t_0), f(t_1)) \leq c$ whenever $t_0, t_1 \in \V(\Lambda)$ are adjacent in $\Lambda$. Set $W = V \cup V^2 \cup \cdots \cup V^c$. Then $W$ is a finite generating set for $G$. Consider $\Cay(G; W)$. Clearly we have that $(f(t_0), f(t_1)) \in \E(\Cay(G; W))$ whenever $(t_0, t_1) \in \E(\Lambda)$. Thus, if we define $\Phi'$ to be the graph with vertex set $G$ and edge set $\{f(e) \: e \in \E(\Lambda)\}$, then $\Phi'$ is a subgraph of $\Cay(G; W)$. Since $f$ is bijective with $\V(\Phi) = G$, $\Phi'$ is graph isomorphic to $\Lambda$ and is a regular spanning tree of $\Cay(G; W)$ of degree $k$.

Now suppose that $G$ has a Cayley graph, $\Cay(G)$, admitting a regular spanning tree $\Phi$ of degree $k$. Then $\Phi$ is non-amenable in the sense that there are $C > 1$ and $m \in \N$ so that $|\rmN_m^\Phi(S)| \geq C \cdot |S|$ for all finite $S \subseteq \V(\Phi)$. Since $\Phi$ is a spanning subgraph of $\Gamma$, it follows that $|\rmN_m^\Gamma(S)| \geq C \cdot |S|$ for all finite $S \subseteq \V(\Gamma)$ (since $\rmN_m^\Gamma(S) \supseteq \rmN_m^\Phi(S)$). Thus, by F{\o}lner's characterization of amenability, $G$ is not amenable.
\end{proof}

The fact that every finitely generated infinite group has some Cayley graph admitting a regular spanning tree is a rather peculiar fact, for as the following proposition shows, not all Cayley graphs have this property.

\begin{prop}
There is a finitely generated infinite group $G$ and a Cayley graph $\Cay(G)$ of $G$ which does not have any regular spanning tree.
\end{prop}

\begin{proof}
Let $G = \Z * \Z_3$. Say $\Z = \langle t \rangle$ and $\Z_3 = \langle u \rangle$. Let $S = \{t, u\}$. Then $S$ is a finite generating set for $G$; in fact it is arguably the most natural generating set for $G$. We claim that $\Cay(G; S)$ does not contain any regular spanning trees.

Due to the structure of free products, it is not difficult to draw the Cayley graph $\Cay(G; S)$. In fact, $\Cay(G; S)$ can be obtained by starting with a regular tree of degree $6$, replacing each of the vertices with a triangle, and evenly dividing the $6$ edges adjacent to the original vertex between the $3$ vertices of the new triangle in its place. Thus every vertex of $\Cay(G; S)$ lies in a triangle (corresponding to repeated multiplication on the right by $u$), and is adjacent to two other triangles (by multiplying on the right by $t$ and $t^{-1}$). If the triangles in $\Cay(G; S)$ are collapsed, then a $6$--regular tree is obtained. A finite portion of $\Cay(G; S)$ is displayed in Figure \ref{Fig Cay}.

\begin{figure}
\begin{center}
\includegraphics[scale=0.2]{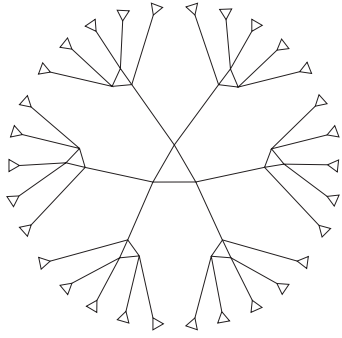}
\caption{A finite portion of the Cayley graph $\Cay(G; S)$}\label{Fig Cay}
\end{center}
\end{figure}

Now let $\Lambda$ be a spanning tree of $\Cay(G; S)$. It suffices to show that $\Lambda$ is not regular. Define
\begin{gather*}
E_t = \{(g, g t) \: g \in G\} \cup \{(g, g t^{-1}) \: g \in G\};\\
E_u = \{(g, g u) \: g \in G\} \cup \{(g, g u^{-1}) \: g \in G\}.
\end{gather*}
Notice that $\E(\Cay(G; S))$ is the disjoint union $E_t \cup E_u$. From the description in the previous paragraph, it is clear that $\E(\Lambda)$ must contain $E_t$ as otherwise $\Lambda$ would either not be spanning or not be connected. Notice that every vertex of $\Cay(G; S)$ is adjacent to precisely two edges from $E_t$. Now consider the three edges $(1_G, u)$, $(u, u^2) = (u, u^{-1})$, and $(u^{-1}, 1_G)$. In order for $\Lambda$ to be spanning, connected, and a tree, $\E(\Lambda)$ must contain precisely two of these three edges. Thus for some $i = -1, 0, 1$ we have $(u^i, u^{i+1}), (u^i, u^{i-1}) \in \E(\Lambda)$ and $(u^{i+1}, u^{i-1}) \not\in \E(\Lambda)$. So $\deg_\Lambda (u^i) = 4$ and $\deg_\Lambda (u^{i+1}) = \deg_\Lambda (u^{i-1}) = 3$. Thus $\Lambda$ is not regular. We conclude that $\Cay(G; S)$ does not contain any regular spanning trees.
\end{proof}

Looking back at the clauses of Corollary \ref{COR SPAN TREE}, we see that the above proposition only shows that ``a Cayley graph'' cannot be replaced by ``every Cayley graph'' in clause (ii). However, it is unclear if clause (i) remains true if ``a Cayley graph'' is replaced by ``every Cayley graph.''

\begin{prob}
Does every infinite Cayley graph with finitely many ends admit a regular spanning tree?
\end{prob}

\begin{prob}
If $G$ is a finitely generated non-amenable group with finitely many ends, then does every Cayley graph of $G$ admit a regular spanning tree of degree strictly greater than two?
\end{prob}

Corollary \ref{COR SPAN TREE} says that if $G$ is finitely generated and non-amenable, then for every $k \geq 3$ $G$ there is a lipschitz bijection from the $k$--regular tree onto $G$. In general these maps cannot be improved to being bilipschitz since some non-amenable groups have finitely many ends and the number of ends is invariant under bilipschitz equivalences. However, it is unclear when we can get a bilipschitz map at the expense of losing bijectivity. By a result of Benjamini--Schramm in \cite{BS}, every non-amenable group contains a bilipschitz image of a regular tree of degree $3$. However, it is not clear if any amenable group contains a bilipschitz image of a regular tree of degree $3$.

\begin{prob}
Which finitely generated groups contain a bilipschitz image of the $3$--regular tree?
\end{prob}

Notice that by Theorem \ref{TWO TREES}, we can equivalently consider bilipschitz images of $k$--regular trees for any $k \geq 3$.

\section{Tilings of Groups} \label{SEC TILE}

In this section we prove that every countable group is poly-MT and every finitely generated group is poly-ccc. We refer the reader to the introduction for the relevant definitions and background.

We first consider poly-ccc groups. As the reader may notice, Corollary \ref{COR TRANS ACT FREE} plays a critical role in these proofs. It is unknown to the author how to prove these results without applying Corollary \ref{COR TRANS ACT FREE}.

\begin{thm} \label{THM POLYCCC}
Every finitely generated group is poly-ccc.
\end{thm}

\begin{proof}
Let $G$ be a finitely generated group. If $G$ is finite then define $T^n = G$ and $\Delta^n = \{1_G\}$ for all $n \in \N$. Then $(\Delta^n; T^n)_{n \in \N}$ is a ccc sequence of monotilings. Thus $G$ is ccc and in particular poly-ccc. Now suppose that $G$ is infinite. Let $S$ be a finite generating set for $G$, and let $d$ be the corresponding left-invariant word length metric. For $m > 1$ let $\F_m$ be the non-abelian free group of rank $m$, and let $\F_1 = \Z$ be the group of integers. By Corollary \ref{COR TRANS ACT FREE} there is a transitive translation-like action $*$ of $\F_m$ on $G$ for some $m \geq 1$.

In \cite{GJS}, Gao, Jackson, and Seward prove that $\F_m$ is a ccc group. Thus there is a ccc sequence of monotilings $(\Delta^n; T^n)_{n \in \N}$ of $\F_m$. Define a map $\phi\co \F_m \rightarrow G$ by setting $\phi(t) = 1_G * t$. Since $\F_m$ acts on $G$ transitively and translation-like, $\phi$ is a bijection. Furthermore, if $\rho$ is the standard left-invariant word length metric on $\F_m$ (associated to the standard generating set of $\F_m$), then $\phi\co (\F_m, \rho) \rightarrow (G, d)$ is lipschitz (Lemma \ref{TRANS ORB LIP}). Say the lipschitz constant is $c$ so that
$$d(\phi(t_1), \phi(t_2)) \leq c \cdot \rho(t_1, t_2)$$
for all $t_1, t_2 \in \F_m$. Now fix $n \in \N$, and consider the family of sets
$$\{\phi(\delta)^{-1} \phi(\delta T^n) \: \delta \in \Delta^n\}.$$
We have that for every $t \in T^n$
$$d(1_G, \phi(\delta)^{-1} \phi(\delta t)) = d(\phi(\delta), \phi(\delta t)) \leq c \cdot \rho(\delta, \delta t) = c \cdot \rho(1_{\F_m}, t).$$
Therefore the family of sets $\{\phi(\delta)^{-1} \phi(\delta T^n) \: \delta \in \Delta^n\}$ all lie within a common ball of finite radius about $1_G \in G$. Since balls of finite radius in $(G, d)$ are finite, the family of sets is finite, say
$$\{\phi(\delta)^{-1} \phi(\delta T^n) \: \delta \in \Delta^n\} = \{T_1^n, T_2^n, \ldots, T_{k(n)}^n\}.$$
Notice that since $\phi$ is injective, all of the $T_i^n$'s have the same number of elements.

Recall that $1_{\F_m} \in \Delta^n$. By reordering the indices if necessary, we may assume that $T_1^n = \phi(1_{\F_m})^{-1} \phi(1_{\F_m} T^n) = \phi(T^n)$. Now define for each $1 \leq i \leq k(n)$
$$\Delta_i^n = \{\phi(\delta) \: \delta \in \Delta^n \text{ and } \phi(\delta)^{-1} \phi(\delta T^n) = T_i^n\}.$$
Notice that $1_G \in \Delta_1^n$. Also notice that if $\delta \in \Delta^n$ and $\phi(\delta) \in \Delta_i^n$, then
$$\phi(\delta) T_i^n = \phi(\delta T^n).$$
The monotiling $(\Delta^n; T^n)$ induces a partition of $\F_m$, and this partition is mapped forward by $\phi$ to a partition of $G$. The above expression shows us that the classes of this partition of $G$ are precisely the left $\Delta_i^n$ translates of $T_i^n$ for $1 \leq i \leq k(n)$. Thus $P_n = (\Delta_1^n, \ldots, \Delta_{k(n)}^n; T_1^n, \ldots, T_{k(n)}^n)$ is a fair polytiling of $G$.

For $n \in \N$ the partition of $\F_m$ induced by $(\Delta^n; T^n)$ is finer than the partition induced by $(\Delta^{n+1}, T^{n+1})$, and thus the image of the first partition is finer than the image of the second. The images of these partitions are precisely the partitions of $G$ induced by $P_n$ and $P_{n+1}$, respectively. So $P_n$ induces a finer partition than $P_{n+1}$. Thus the sequence $(P_n)_{n \in \N}$ is coherent. Also, as previously mentioned, $1_G \in \Delta_1^n$ and $T_1^n = \phi(T^n) \subseteq \phi(T^{n+1}) = T_1^{n+1}$ for all $n \in \N$. Since $\phi$ is bijective and $\F_m = \bigcup T^n$, we have that $G = \bigcup T_1^n$. Thus $(P_n)_{n \in \N}$ is a coherent, centered, and cofinal sequence of fair polytilings. We conclude that $G$ is poly-ccc.
\end{proof}

The above theorem shows that a very large class of groups are poly-ccc. However, we don't know if all countable groups are poly-ccc. We know even less about which groups are ccc.

\begin{prob}
Is every countable group poly-ccc? Is every countable group ccc?
\end{prob}

\begin{cor} \label{COR PMT}
Every countable group is poly-MT.
\end{cor}

\begin{proof}
Let $G$ be a countable group, and let $F \subseteq G$ be finite. Then $H = \langle F \rangle$ is a finitely generated group and is thus poly-ccc by the previous theorem. Thus there is a fair sequence of ccc polytilings
$$(P_n)_{n \in \N} = (\Delta_1^n, \ldots, \Delta_{k(n)}^n; T_1^n, \ldots, T_{k(n)}^n)_{n \in \N}$$
of $H$. By definition, $T_1^n \subseteq T_1^{n+1}$ and $H = \bigcup T_1^n$. Thus there is $n \in \N$ with $F \subseteq T_1^n$. Fix this value of $n$. Let $D$ be a complete set of representatives for the left cosets of $H$ in $G$. Specifically, $D H = G$ and $d H \cap d' H = \varnothing$ for $d \neq d' \in D$. It is easy to see that $(D \Delta_1^n, \ldots, D \Delta_{k(n)}^n; T_1^n, \ldots, T_{k(n)}^n)$ is a fair polytiling of $G$. Since $F \subseteq T_1^n$, we conclude that $G$ is poly-MT.
\end{proof}

\begin{prob}
Is every countable group MT?
\end{prob}

Below we address what values $|T_1| = |T_2| = \cdots = |T_k|$ one can get from fair polytiles $(T_1, \ldots, T_k)$. Recall that a group is \emph{locally finite} if every finite subset generates a finite subgroup. If $G$ is a locally finite group and $(T_1, \ldots, T_k)$ is a fair polytile of $G$, then $|T_1|$ must divide the order of $H = \langle \bigcup_{i = 1}^k T_i \rangle$. This is because any set $a T_i \subseteq a H$ must either be contained in $H$ or be disjoint from $H$. So $(T_1, \ldots, T_k)$ is a fair polytile for the finite group $H$, and thus $|T_1|$ divides $|H|$. Thus we see that in the case of locally finite groups there are restrictions on what the cardinality of the sets $T_i$ can be. The following theorem shows that this is the only obstruction.

\begin{thm}
Let $G$ be a countable non-locally finite group. Then for every $n \geq 1$ there is a fair polytile $(T_1, \ldots, T_k)$ of $G$ with $|T_1| = n$. Furthermore, if $F \subseteq G$ is finite then for all sufficiently large $n$ there is a fair polytile $(T_1, \ldots, T_k)$ of $G$ with $F \subseteq T_1$ and $|T_1| = n$.
\end{thm}

\begin{proof}
Let $F \subseteq G$ be finite. Since $G$ is not locally finite, there must exist a finite set $S \subseteq G$ with $\langle S \rangle$ infinite. Set $H = \langle F \cup S \rangle$. So $H$ is a finitely generated infinite group. By Corollary \ref{COR TRANS ACT FREE} there is a transitive translation-like action $*$ of $\F_m$ on $H$ for some $m \geq 1$. Define $\phi\co \F_m \rightarrow H$ by $\phi(t) = 1_H * t$, and let $\Cay(\F_m)$ be the canonical Cayley graph of $\F_m$. In \cite{GJS}, it is proven that if $\Phi$ is a finite connected subgraph of $\Cay(\F_m)$ then $\V(\Phi)$ is a monotile of $\F_m$. In particular, for every $n \geq 1$ there is a monotile $T$ of $\F_m$ with $|T| = n$. It also follows that for every finite set $F' \subseteq \F_m$ and for all sufficiently large $n$ there is a monotile $T$ of $\F_m$ with $F' \subseteq T$ and $|T| = n$. The statement of this theorem follows by using $F' = \phi^{-1}(F)$, then applying the construction appearing in the proof of Theorem \ref{THM POLYCCC}, and finally using left coset representatives for left cosets of $H$ in $G$ as in the proof of Corollary \ref{COR PMT}.
\end{proof}

We have shown that all countable groups are poly-MT and all finitely generated groups are poly-ccc, so it seems quite possible that all countable groups are poly-ccc. However, there are other natural tiling properties one can consider which lie between ccc and poly-ccc and lie between MT and poly-MT. We do not study these notions here, but we define them below so that others may investigate them.

\begin{defn}
A countable group $G$ is \emph{super poly-MT} if for every finite $F \subseteq G$ there is a fair polytile $(T_1, \ldots, T_k)$ satisfying $F \subseteq \bigcap T_i$.
\end{defn}

\begin{defn}
A countable group $G$ is \emph{super poly-ccc} if $G$ has a fair ccc sequence of polytilings
$$(\Delta_1^n, \ldots, \Delta_{k(n)}^n; T_1^n, \ldots, T_{k(n)}^n)$$
such that for each $n$ we have $\bigcap T_i^n \subseteq \bigcap T_i^{n+1}$ and $G = \bigcup_n \bigcap_i T_i^n$.
\end{defn}

The following diagram clarifies the relationship between the six tiling properties discussed in this paper.

$$\begin{array}{ccc}
\text{ccc} & \Longrightarrow & \text{MT} \\
\Downarrow & & \Downarrow \\
\text{super poly-ccc} & \Longrightarrow & \text{super poly-MT} \\
\Downarrow & & \Downarrow \\
\text{poly-ccc} & \Longrightarrow & \text{poly-MT}
\end{array}$$

The ccc property is the strongest of all of the mentioned tiling properties, and there are currently no groups which are known to not be ccc. However, showing that groups are ccc can be quite difficult, and so studying some of the other related notions in the diagram may be more fruitful.

\begin{prob}
Which groups are super poly-ccc?
\end{prob}

\begin{prob}
Which groups are super poly-MT?
\end{prob}

\end{document}